\documentclass{article}

\usepackage{PRIMEarxiv}

\usepackage[utf8]{inputenc} 
\usepackage[T1]{fontenc}    
\usepackage{hyperref}       
\usepackage{url}            
\usepackage{booktabs}       
\usepackage{amsfonts}       
\usepackage{nicefrac}       
\usepackage{microtype}      
\usepackage{lipsum}
\usepackage{fancyhdr}       
\usepackage{graphicx}       
\graphicspath{{media/}}     

\title{A Lattice Boltzmann Method for Elastic Solids Under Plane Strain Deformation
}

\author{
  Alexander Schl\"{u}ter, Sikang Yan, Erik Faust \\
  Institute of Applied Mechanics \\
  Technische Universit\"{a}t Kaiserslautern\\
  D-67653, Kaiserslautern\\
  \texttt{\{Alexander Schl\"{u}ter\} aschluet@rhrk.uni-kl.de} \\
   \And
  Henning M\"{u}ller, Ralf M\"{u}ller \\
 Institut f\"{u}r Mechanik \\
  Technische Universität Darmstadt\\
  D-64287, Darmstadt\\
}

\renewcommand{\vec}{\boldsymbol}        
\newcommand{\eps}{\varepsilon}
\newcommand{\cc}{\boldsymbol{c}}
\newcommand{\argum}{(\vec{x},t)}
\newcommand{\norm}[1]{\left\lVert#1\right\rVert}

\usepackage{nicefrac}
\usepackage[]{graphicx}	
\usepackage{amsmath}
\graphicspath{ {pics/} }
\usepackage{xcolor}

\begin{document}
\maketitle

\begin{abstract}
The Lattice Boltzmann Method (LBM), e.g. in \cite{chopard_lattice_1999} and \cite{kruger_lattice_2017}, can be interpreted as an alternative method for the numerical solution of partial differential equations. Consequently, although the LBM is usually applied to solve fluid flows, the above interpretation of the LBM as a general numerical tool, allows the LBM to be extended to solid mechanics as well. In this spirit, the LBM has been studied in recent years. First publications \cite{schluter_lattice_2018}, \cite{reinirkens_lattice_2018} presented an LBM scheme for the numerical solution of the dynamic behavior of a linear elastic solid under simplified deformation assumptions. For so-called anti-plane shear deformation, the only non-zero displacement component is governed by a two-dimensional wave equation. In this work, an existing LBM for the two-dimensional wave equation is extended to more general plane strain problems. The proposed algorithm reduces the plane strain problem to the solution of two separate wave equations for the volume dilatation and the non-zero component of the rotation vector, respectively. A particular focus is on the implementation of types of boundary conditions that are commonly encountered in engineering practice for solids: Dirichlet and Neumann boundary conditions. Last, several numerical experiments are conducted that highlight the performance of the new LBM in comparison to the Finite Element Method.
\end{abstract}


\section{Introduction}
\label{sec:intro}

The mechanical behavior of solid bodies is of interest to both engineering and science. Thus, a large number of numerical methods capable of dealing with elasticity have emerged over time. The more prominent ones among these, finite differences methods (FDM), finite element methods (FEM) and finite volume methods (FVM), work on the principle of discretizing the domain of interest and replacing the governing system of differential equations by algebraic equations. Such methods take a kind of top-down approach, and can therefore be thought of as acting on a macroscopic scale.
In contrast, some numerical methods, such as molecular dynamics (MD) or density functional theory (DFT), regard the interactions of a system's most basic constituents, such as individual particles and electrons, on a microscopic scale.

A different approach is taken with \emph{Lattice-Boltzmann methods} (LBMs). The common principle of this type of methods is to transform the given physical problem into a transport problem. Based on Boltzmann's transport equation from statistical mechanics, distribution functions are transported across phase-space, which is discretized both by a regular lattice and a set of associated lattice velocities.
Information is exchanged between neighboring lattice sites in a streaming-like process along links connecting these points.
This information is represented by so-called distribution functions, where each distribution function is associated with a different lattice velocity. The distribution functions are subjected to on-site interactions, or \textsl{collisions}, which incorporate the underlying microscopic theory in a probabilistic manner. Thus, LBMs can be said to act mesoscopically, i.e. on an intermediate scale.

LBMs are well established in computational fluid dynamics (CFD) and have subsequently been extended to further scientific fields, such as solving Schrödinger's equation \cite{succi_numerical_1996} or Wigner's equation \cite{solorzano_lattice_2018} in quantum mechanics.
Developing a LBM for solid mechanics could mean using a single method on both sides of a fluid-solid-interface, which is a topic of interest \cite{bungartz_fluid-structure_2006} in CFD.

We approach the topic of a LBM for elastic bodies from the mechanical point of view. This includes a greater focus on finite domains with an appropriate boundary handling, which is of great concern for engineering problems. 
The advantages over the established methods in computational engineering include the generally great computational efficiency, while still being able to handle the boundary conditions of complex domains. A further improvement of computation times by can be achieved by employing parallel computing, which is easy to implement with most LBMs. This opens up possibilities in highly dynamical problems, requiring very fine resolution of the temporal domain.

The dynamic behavior of  elastic solids can be described by multiple wave equations, that are superposed to obtain the aggregated deformation. This mathematical description is closer to the transport phenomena for which the LBM was initially conceived, when compared to the original Navier-Cauchy equation. In fact, LBMs have already been developed for the wave equation, see e.g.~\cite{yan_lattice_2000,chopard_cellular_1998,frantziskonis_lattice_2011}.
Furthermore, LBMs have been proposed for the numerical treatment of mechanical problems in solid bodies~\cite{chopard_lattice_1998,chopard_lattice_1999,marconi_lattice_2003,mora_lattice_1992,yin_direct_2016,murthy_lattice_2018,escande_lattice_2020}, and specifically elastic wave propagation~\cite{dhuri_numerical_2017,jiang_acoustic_2019}, which is also a topic of interest in geophysics and seismology. However, an extensive method for the deformation of linear elastic solids under loads, containing appropriate boundary conditions, has still to be accomplished.

In previous works~\cite{schluter_lattice_2018,schluter_boundary_2021} we applied the LBM for wave equations published by Yan~\cite{yan_lattice_2000} to the mechanical problem of anti-plane shear deformation, which we then used for fracture mechanics. This work now regards the two-dimensional problem of plane strain.
The fundamental idea is to decompose the plane strain problem governed by the Navier-Cauchy equation into two equivalent wave equations that are solved with the LBM
for wave propagation  by Chopard et al.~\cite{chopard_lattice_1998}.

The discussion is structured as follows:
First the mechanical problem is reviewed and the relevant equations are derived.
The next section introduces the LBM for the wave equation, followed by the presentation of the algorithm for the plane strain case. This includes a treatment of boundary conditions similar to~\cite{schluter_boundary_2021}.
Lastly three numerical examples show the feasibility of our algorithm, each compared to FEM computations, which act as benchmarks.

\section{Plane Strain Deformation of a Linear Elastic Solid}
\label{sec:elastic}

\begin{figure}[bt]
    \centering
    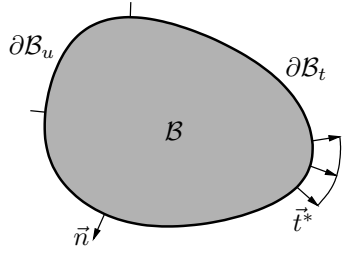
    \caption{%
        Body $\mathcal{B}$ with outer normal vector $\boldsymbol{n}$, subjected to Neumann boundary conditions $\vec{t}^*$ on $\mathcal{B}_t$ and Dirichlet boundary conditions on $\mathcal{B}_u$.
    }
    \label{fig:body}
\end{figure}

We consider a homogeneous, isotropic, and elastic body $\mathcal{B}$ with boundary $\partial\mathcal{B}=\partial\mathcal{B}_u\cup\mathcal{B}_t$, which is subjected to Dirichlet boundary conditions $\boldsymbol{u}=\vec{u}^*$ for the displacement $\vec{u}$  on $\partial\mathcal{B}_u$ and Neumann boundary conditions $\vec{\sigma}\vec{n}=\boldsymbol{t}^*$ for the Cauchy stress tensor $\vec{\sigma}$ on $\partial\mathcal{B}_t$, see  Fig.~\ref{fig:body}. For the plane strain case the problem is only regarded in two dimensions and under small strain assumptions.

A set of fundamental equation is taken as a basis for the derivation of the mathematical description of the problem.
Firstly, the strain-displacement relation for small strains is given by the linearized strain tensor
\begin{equation}
    \vec{\eps} = \frac{1}{2} \left( \nabla \vec{u} + (\nabla \vec{u})^{T} \right),  \label{eq:strain}
\end{equation}
where $\vec{u} = \vec{u}(x,y,t)$ describes the time-dependent displacement field in two dimensions under plane strain assumptions.
The general equation of motion for the small strain case, in the absence of a body force, is given by
\begin{equation}
    \nabla \cdot \vec{\sigma} = \rho \, \frac{\partial^{2} \vec{u}}{\partial t^{2}},    \label{eq:cauchy}
\end{equation}
with Hooke's law as the linear stress-strain relation
\begin{equation}
    \vec{\sigma} = \lambda \, \mathrm{tr}(\vec{\eps}) \, \boldsymbol{1} + 2 \mu \vec{\eps}. \label{eq:hooke}
\end{equation}
Herein $\lambda$ and $\mu$ are the Lam\'{e} parameters of the material, $\boldsymbol{1}$ is the second-order identity tensor and the operator $\mathrm{tr}(*)$ denotes the trace of a second order tensor.

Equation~\eqref{eq:strain} is substituted in equation~\eqref{eq:hooke}, which is then substituted in~\eqref{eq:cauchy}. The result is the Navier-Cauchy equation
\begin{equation}
    (\lambda + \mu) \, \nabla \left( \nabla \cdot \vec{u} \right) + \mu \nabla^2 \, \vec{u} = \rho \, \frac{\partial^{2} \vec{u}}{\partial t^{2}},    \label{eq:navier1}
\end{equation}
which describes the mechanical behavior of an isotropic linear elastic solid.
Using the general identity 
\begin{equation}
    \nabla^2 \vec{u} = \nabla (\nabla \cdot \vec{u}) - \nabla \times (\nabla \times \vec{u})
\end{equation}
from vector calculus, 
equation~\eqref{eq:navier1} can be rewritten as 
\begin{equation}
    c_d^2 \, \nabla (\nabla \cdot \vec{u}) - c_s^2 \, \nabla \times (\nabla \times \vec{u}) = \frac{\partial^{2} \vec{u}}{\partial t^{2}},   \label{eq:wave-navier}
\end{equation}
where $c_d = \sqrt{\nicefrac{(\lambda + 2 \mu)}{\rho}}$ and $c_s = \sqrt{\nicefrac{\lambda}{\rho}}$.

With regard to equation~\eqref{eq:wave-navier}, two fields, $\phi$ and $\vec{\psi}$, can be defined as follows:
\begin{align}
    \phi = \nabla \cdot \vec{u} &&\text{and} &&\vec{\psi} = \nabla \times \vec{u}. \label{eq:phi-psi}
\end{align}
The scalar field $\phi$ describes the dilatation of the displacement field $\vec{u}$, whereas the vector field $\vec{\psi}$ describes the rotation of $\vec{u}$. In two dimensions, the latter reduces to $\vec{\psi} = \psi \, \vec{e}_z$.
The Navier-Cauchy equation can be restated in terms of these fields
\begin{equation}
    c_d^2 \, \nabla \phi - c_s^2 \, \nabla \times \vec{\psi} = \frac{\partial^{2} \vec{u}}{\partial t^{2}}.   \label{eq:navier} 
\end{equation}

In conjunction with the definitions in equation~\eqref{eq:phi-psi}, applying the divergence to both sides of equation~\eqref{eq:navier} results in
\begin{subequations}
    \label{eq:wave-eqs}
    \begin{equation}
        c_d^2 \, \nabla^2 \phi = \frac{\partial^{2} \phi}{\partial t^{2}} \label{eq:dil-wave},
    \end{equation}
    while applying the curl results in 
    \begin{equation}
        c_s^2 \, \nabla^2 \vec{\psi} = \frac{\partial^{2} \vec{\psi}}{\partial t^{2}}.  \label{eq:rot-wave}
    \end{equation}
\end{subequations}
Thus, the Navier-Cauchy equation~\eqref{eq:navier} can be reduced to two wave equations, the dilatational wave equation~\eqref{eq:dil-wave} for ${ \phi = \nabla \cdot \vec{u} }$ with wave speed $c_d$, and the rotational wave equation~\eqref{eq:rot-wave} for ${ \vec{\psi} = \psi \, \vec{e}_z = \nabla \times \vec{u} }$ with wave speed $c_s$.
Note that there are other ‘decompositions' of the displacement field besides~\eqref{eq:phi-psi} that lead to similar wave equations, see e.g.~\cite{sternberg_integration_1960}.

\section{Lattice Boltzmann Method for Plane Strain}
\label{sec:planeStrainlbm}

\begin{figure}[bt]
    \centering
    \vspace{0.1cm}
    \input{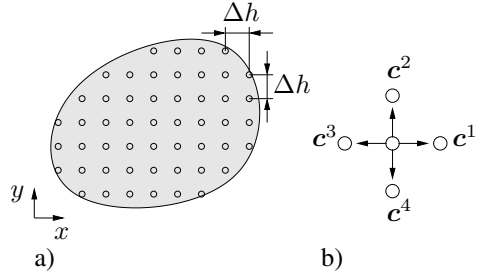}
    \caption{a) Lattice representation of the elastic solid and b)  and the associated lattice velocity vectors (lattice links) for a single lattice point.}
\label{fig:domain}
\end{figure}

The proposed numerical strategy for the plane strain case relies on solving the wave equations~(\ref{eq:wave-eqs}) by means of the LBM by Chopard et al. \cite{chopard_lattice_1998}. In the LBM,  a body $\mathcal{B}$ is typically approximated by a regular lattice with lattice spacing $\Delta{h}$ as depicted in Fig.~\ref{fig:domain}~a). The approach by Chopard et al. is based on a D2Q5 scheme, see also Fig.~\ref{fig:domain}, with the lattice velocities 
\begin{gather}
    \begin{aligned}
        \vec{c}^{0} &=({c}^{0}_x,{c}^{0}_y) = (0,0) \\
        \vec{c}^{1} &=({c}^{1}_x,{c}^{1}_y) = (c,0) \\
        \vec{c}^{2} &=({c}^{2}_x,{c}^{2}_y) = (0,c) \\
        \vec{c}^{3} &=({c}^{3}_x,{c}^{3}_y) = (\text{-}c,0) \\
        \vec{c}^{4} &=({c}^{4}_x,{c}^{4}_y) = (0,\text{-}c)
    \end{aligned}
\end{gather}
where $c=\nicefrac{\Delta{h}}{\Delta{t}}$ is the speed at which information can travel in the lattice. Thus, the lattice velocities $\vec{c}^{1}, \vec{c}^{2}, \vec{c}^{3}, \vec{c}^{4}$ allow for information to be transported  to each of the four neighbors of a lattice point in a so-called D2Q5 scheme\footnote{D2Q5 refers to the dimension of the lattice, i.e. two in this case, and the number of lattice velocities, i.e. five in this case.} in one time step, whereas $\vec{c}^{0}$ is associated with information remaining at a particular lattice point. Information is represented by distribution functions, e.g. $f^{\alpha}$ represents information, which is transported with lattice velocity $\vec{c}^\alpha$.

In order to simulate the wave equations (\ref{eq:wave-eqs}), the distribution functions need to be interpreted, i.e. a relation to the macroscopic fields needs to be established.
We introduce two sets of distribution functions to simulate both wave equations and relate them to the macroscopic fields through
\begin{equation}
   \sum_{\alpha=0}^4 f_{\psi}^{\alpha} = \psi, \quad  \sum_{\alpha=0}^4 f_{\phi}^{\alpha} = \phi. 
   \label{eq:interpretation}
\end{equation}

The Lattice Boltzmann equation (LBE) models transport as well as the interaction of distribution functions between and at lattice points respectively. Since two wave equations need to be solved, we also introduce two associated LBEs for $\psi$ and $\phi$ respectively
\begin{align}
    f^{\alpha}_{\psi\vert\phi} \left( \vec{x} + \vec{c}^{\alpha} \Delta{t}, t + \Delta{t} \right) & =  \nonumber\\
    f^{\alpha}_{\psi\vert\phi} \argum - \frac{\Delta{t}}{\tau} \left[ f^{\alpha}_{\psi\vert\phi} \argum \right.& \left. -   f^{\alpha}_{\text{eq},\psi\vert\phi} \argum \right],
    \label{eq:LBE}
\end{align}
where the notation $(\psi\vert\phi)$ indicates that either $\psi$ or $\phi$ need to be chosen for the whole equation and the common BGK approximation is employed, see~\cite{bhatnagar_model_1954} and~\cite{welander_temperature_1954}.
Equation (\ref{eq:LBE}) is universal to many Lattice Boltzmann models. The specific physics can be modeled by choosing the equilibrium distribution functions $f^{\alpha}_{\text{eq},\psi\vert\phi}$ and relaxation time $\tau$ in a certain way. In order to model a wave equation Chopard et al. propose $\tau=0.5\Delta{t}$ and 
\begin{align}
    &f^{0}_{\text{eq},\psi\vert\phi} = a_{0,\psi\vert\phi} (\psi\vert\phi) \nonumber\\ 
    &f^{\alpha}_{\text{eq},\psi\vert\phi} = a_{\psi\vert\phi}(\psi\vert\phi) + b \dfrac{\vec{c}^\alpha\cdot\vec{J}_{\psi\vert\phi}}{2c^2}, \nonumber\\
    \text{with}\ &\vec{J}_{\psi\vert\phi} = \sum_{\alpha=0}^{4} \vec{c}^{\alpha} f^{\alpha}_{\psi\vert\phi} \argum, \ \text{for }\alpha \neq 0,
    \label{eq:feq}
\end{align}
where again the notation $(\psi\vert\phi)$ indicates that either $\psi$ or $\phi$ need to be chosen for the whole equation.
The parameters also need to fulfill the requirements
\begin{align}
\begin{split}
   & b=1, \quad \text{conservation of } \vec{J}_{\psi\vert\phi}\\
   & a_{0,\psi\vert\phi}+ 4a_{\psi\vert\phi} = 1, \quad \text{conservation of } {\psi\vert\phi}\\
   & a_{0,\psi\vert\phi}\ge{0}, \quad {\text{stability}}
   \end{split}
   \label{eq:feq_req}
\end{align}
and 
\begin{equation}
    c_s\vert{c_d} =\dfrac{\Delta{h}}{\Delta{t}}\sqrt{2a_{\psi\vert\phi}},
    \label{eq:parameters_to_wavespeed_relation}
\end{equation}
see Chopard~\cite{chopard_lattice_1999}.
Equation (\ref{eq:parameters_to_wavespeed_relation}) allows us to adjust the macroscopic wave speed modeled by the LBM independently of the time step $\Delta{t}$ or the lattice spacing $\Delta{h}$ by choosing $a_{\psi\vert\phi}$ accordingly. We exploit this in order to be able to simulate both wave equations (\ref{eq:wave-eqs}) on the same lattice, i.e. fixed $\Delta{h}$, and the same time discretization, i.e. fixed $\Delta{t}$. Apart from (\ref{eq:parameters_to_wavespeed_relation}), the requirements (\ref{eq:feq_req}) still need to be fulfilled.
This can be accomplished for the general case $c_s<c_d$ by setting 
\begin{align}
\begin{split}
0\le{a_{\phi}}\le0.25, \\
a_{\psi} = \dfrac{c_s^2}{c_d^2}a_{\phi},\\
\Delta{t} = \dfrac{\Delta{h}}{c_s}\sqrt{2a_\psi}= \dfrac{\Delta{h}}{c_d}\sqrt{2a_\phi}.
\end{split}
\label{eq:conditions-parameters}
\end{align}
Note that the Courant-Friedrichs-Lewy (CFL) stability condition~\cite{Courant_1928}, which is critical for the numerical analysis of hyperbolic PDEs with explicit schemes, for the larger -- and more critical -- wave speed
\begin{equation}
    \dfrac{c_d\Delta{t}}{\Delta{h}} \le{1},
\end{equation}
is always guaranteed by (\ref{eq:conditions-parameters}).

	\begin{figure*}
\fbox{\parbox{\textwidth}{
\begin{itemize}
\item Preprocessing
\begin{itemize}
    \item Build lattice for geometry
    \item Compute surface measure and cell volumes for boundary lattice points
\end{itemize}
\item Solver
\begin{itemize}
    \item Initialize ${\vec{u}}(\vec{x},0)$ and $\dot{\vec{u}}(\vec{x},0)$
    \item Initialize $\vec{\psi}(\vec{x},0) = \nabla\times\vec{u}(\vec{x},0)$, $\vec{\phi}(\vec{x},0) = \nabla\times\vec{u}(\vec{x},0)$ computed by finite differences
    \item Initialize the distribution functions by (\ref{eq:feq}) 
    \item Start time loop
    \begin{itemize}
        \item ${t}\to{{t}+\Delta{t}}$ 
        \item Compute accelerations in the interior $\ddot{\vec{u}}(\vec{x},t) = c_d^2\nabla{{\phi} \argum} -c_s^2\vec{\psi}\argum$ i.e. by (\ref{eq:wave-navier})
        \item Compute accelerations at boundary points by (\ref{eq:acc_neumann}) or (\ref{eq:acc_dirichlet})
        \item Compute $\vec{u}(\vec{x},t+\Delta{t})$ by explicit integration, i.e. by (\ref{eq:integration})
        \item set $\boldsymbol{\psi}(\vec{x},t+\Delta{t})$ at boundary lattice points to be consistent with $\vec{u}(\vec{x},t+\Delta{t})$ according to (\ref{eq:psi_phi_at_bound})
        \item Solve the wave equations for ${\psi}(\vec{x},t+\Delta{t})$ and $\phi(\vec{x},t+\Delta{t})$ via the LBE (\ref{eq:LBE}) and the interpretation
        (\ref{eq:interpretation})
        \item Every $l$-th time step perform synchronization, see section~\ref{sec:synchronization}
    \end{itemize}
    \item End time loop if ${t}={t}_{\text{final}}$
\end{itemize}
\end{itemize}}}
\caption{Summary of the employed lattice Boltzmann algorithm.}
\label{fig:algorithm}
\end{figure*}

The LBE (\ref{eq:LBE}) is only part of the overall algorithm that is used to solve a plane strain problem  as summarized in Fig.~\ref{fig:algorithm}. The other parts of the algorithm are discussed in the following.

\begin{figure}
    \centering
    \includegraphics[width=0.4\textwidth]{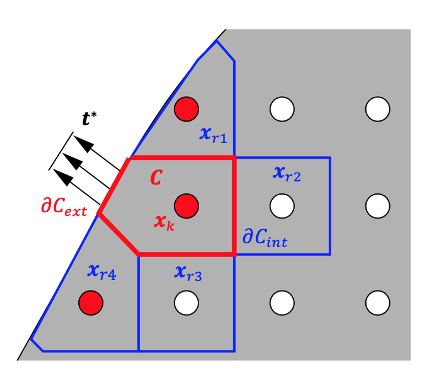}
    \caption{Cells are generated around each boundary lattice point in order to apply Neumann boundary conditions.}
    \label{fig:boundary_neumann}
\end{figure}

The initial preprocessing step builds the lattice for a given geometry and also computes cells at each boundary lattice point as depicted in Fig.~\ref{fig:boundary_neumann}. 
Without going into detail, the algorithm for creating individual cells starts with a quadratic cell of side length $\Delta{h}$ centered around a boundary lattice point. This original cell is subsequently modified to match the boundary geometry.
The volume $V_C$ of a cell and the surface that such a cell $C$ shares with the external boundary $\partial{C}_{\text{ext}}\subset{\partial{\cal B}}$ are relevant for the computation of the acceleration at boundary lattice points later on.

After preprocessing, the material velocity $\dot{\boldsymbol{u}}$ and the displacement $\vec{u}$ are initialized first. Subsequently, $\psi$ and $\phi$ are initialized by a finite difference approximation of (\ref{eq:phi-psi}). Lastly,  the initial distribution functions are determined to be the value of the equilibrium distribution function
\begin{equation}
    f^{\alpha}(\vec{x},0)_{\psi\vert\phi} = f^{\alpha}_{\text{eq},\psi\vert\phi}(\psi(\vec{x},0)\vert\phi(\vec{x},0)).
\end{equation}

In the time loop, the acceleration is computed from the Navier-Cauchy equation (\ref{eq:navier}) at interior lattice points, whereas boundary conditions determine the acceleration  at the boundary lattice points. Once the acceleration at each lattice point is known, the displacement is computed by explicit integration via the Newmark method, i.e.
\begin{align}
    \vec{u}(\vec{x}, t+\Delta{t}) = \vec{u}(\vec{x}, t) + \Delta{t}\dot{\vec{u}}\argum + \dfrac{\Delta{t}^2}{2} \ddot{\vec{u}}(\vec{x},t) \nonumber,\\
    \dot{\vec{u}}(\vec{x}, t+\Delta{t}) = \dot{\vec{u}}(\vec{x}, t) + \Delta{t}\,\ddot{\vec{u}}\argum. 
    \label{eq:integration}
\end{align}
After the integration step, the displacement field is already updated, i.e. $\vec{u}(\vec{x}, t+\Delta{t})$ is determined at all lattice points. The rotation $\psi\argum$, the dilatation $\phi\argum$ and the associated distribution functions $f_{\psi\vert\phi}^\alpha$ have not been updated yet, but they are required to compute the acceleration at interior lattice points in the next time step.

We prepare the required update by computing rotation and dilatation fields as well as distribution functions at the boundary lattice points. All of these must be consistent with the applied boundary conditions as well, see the next section for details. 
Subsequently, the rotation and dilatation fields are updated in the interior by the LBM. This step includes the update of all interior distribution functions by (\ref{eq:LBE}) and the computation of the rotation and the dilatation by (\ref{eq:interpretation}).

\subsection{Boundary Conditions}

In this section, the treatment of boundary conditions is explained in more detail since it is the most complex part of the proposed LBM. The overall strategy is to first determine the acceleration $\ddot{\vec{u}}\argum$, that is consistent with the boundary conditions for each boundary lattice point. Subsequently, the displacement field is updated everywhere via the Newmark integration mentioned above. Last, the rotation and dilatation fields as well as the distribution functions are reconciled with the displacement. 

\subsubsection{Consistent Acceleration at Boundary Lattice Points}
For Neumann type boundary conditions the consistent acceleration is determined by computing cells $C$ with size $V_C$ and boundary $\partial{C}$ around each boundary lattice point $\vec{x}_k$ as shown in Fig.~\ref{fig:boundary_neumann}. For each of these cells, we consider a balance of momentum
\begin{align}
\int_{C}\rho \ddot{\vec{u}}(\vec{x}_k,t)\, {\rm d}v & = \nonumber\\
\int_{\partial{C}_{\text{int}}} & \vec{\sigma}(\vec{x},t)\vec{n}\, \text{d}a + \int_{\partial{C}_{\text{ext}}} \vec{t}^*(\vec{x},t)\, \text{d}a,
    \label{eq:contBalanceOfMomentum}
\end{align}
where $\partial{C}_{\text{int}}$ is the part of the boundary of the cell which is shared with neighboring cells and $\partial{C}_{\text{ext}}$ is part of the boundary of the cell that is shared with the boundary of the body.
Equation~(\ref{eq:contBalanceOfMomentum}) is 
simplified by assuming that $\rho$ and $\ddot{\vec{u}}(\vec{x}_k,t)$ are constant across the cell and that the stress $\vec{\sigma}_{kr}$ is constant for each segment of the internal boundary  shared with a particular neighbor $\vec{x}_r$,
\begin{align}
    \ddot{\vec{u}}(\vec{x}_k,t)& \approx \nonumber\\* \dfrac{1}{\rho{V_C}}&\left(\sum_{r\in\text{Neighbors}}\vec{\sigma}_{kr}\vec{n}_{kr} + \int_{\partial{C}_{\text{ext}}}\vec{t}^*(\vec{x},t)\,\text{d}a\right).
    \label{eq:acc_neumann}
\end{align}
The surface measure, i.e. the length of of the boundary segment in 2D and the normal vector are denoted by $l_{kr}$ and $\vec{n}_{kr}$ respectively. The stress tensor at each segment is approximated by 
\begin{equation}
    \vec{\sigma}_{kr} = \dfrac{1}{2}\left(\vec{\sigma}(\vec{x}_k,t)+\vec{\sigma}(\vec{x}_r,t)\right),
\end{equation}
where the stress at the lattice points $\vec{x}_k$ and $\vec{x}_r$ is computed by a finite difference approximation of (\ref{eq:hooke}).

For Dirichlet type boundary conditions $\vec{u}=\vec{u}^*$ on $\partial{B}_{{u}}$, where $\vec{u}^*$ is the prescribed displacement value, the acceleration $\ddot{\vec{u}}(\vec{x}_k,t)$ at boundary lattice points is determined from the integration scheme (\ref{eq:integration}). Although extrapolation to a non-lattice conforming boundary is also possible, we limit the discussion of Dirichlet boundary conditions to situations in which boundary lattice points lie exactly on the boundary. In this case, it is 
\begin{equation}
    \vec{u}(\vec{x}_k,t+\Delta{t}) = \vec{u}^*.
\end{equation}
Thus, (\ref{eq:integration}) can be solved for the required acceleration 
\begin{equation}
    \ddot{\vec{u}}({\vec{x}_k,t}) = \dfrac{2}{\Delta{t}^2}(\vec{u}^*(t)-\vec{u}(\vec{x}_k,t))-\dfrac{2}{\Delta{t}}\dot{\vec{u}}(\vec{x}_k,t).
    \label{eq:acc_dirichlet}
\end{equation}

\subsubsection{Consistent Displacement, Rotation, Dilatation and Distribution Functions at Boundary Lattice Points}

\begin{figure*}[htb]
    \centering
    \includegraphics[width=\textwidth]{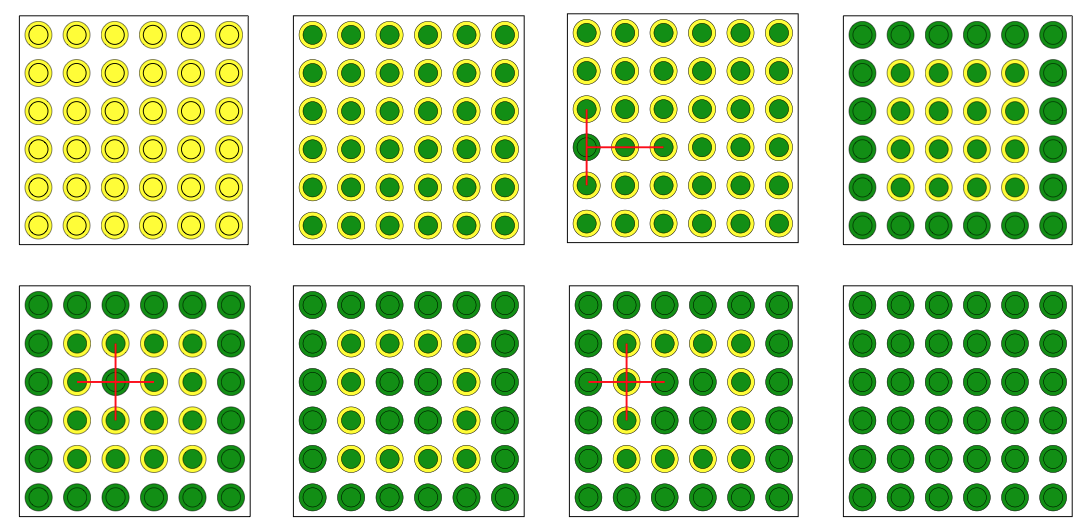}
    \begin{picture}(0,0)
        \put(-220,121){\small a)}
        \put(-105,121){\small b)}
        \put(13,121){\small c)}
        \put(127,121){\small d)}
         \put(-220,5){\small e)}
         \put(-105,5){\small f)}
        \put(13,5){\small g)}
        \put(127,5){\small h)}
    \end{picture}
    \caption{Updating the displacement $\boldsymbol{u}$ as well as the rotation $\psi$ and dilation $\phi$ in a simple square domain. The outer circles represent the state of $\psi$ and $\phi$ at the particular lattice points, whereas the inner circles represent the state of $\boldsymbol{u}$. Yellow indicates that quantities still have a value associated with the previous time step $t$, whereas green color indicates that  $\boldsymbol{u}$ or $\psi$ and $\phi$ are already updated to their values at $t+\Delta{t}$. a) the state after the previous time step. b) integration is performed at all lattice points which updates the displacement. c) finite differences (stencil is indicated by the red lines) are used to update $\psi$ and $\phi$ at the boundary lattice points in a way that is consistent with the new displacement field. d) all boundary points are updated. e) the LBM update, i.e. solving the wave equations,  leads to a consistent rotation and dilatation at interior lattice points (red lines indicate from which neighbors information is streamed to an interior lattice point). f) all interior points have consistent fields after the LBM update. g) intermediate `second row' boundary points are more problematic since there is also information streamed from boundary points that does not originate from the LBE for the wave equations, but from the handling of boundary conditions (red lines indicate from which neighbors information is streamed to an interior lattice point). h) Fields are consistent -- considering the previous remarks -- at all lattice points.}
    \label{fig:boundary_update}
\end{figure*}
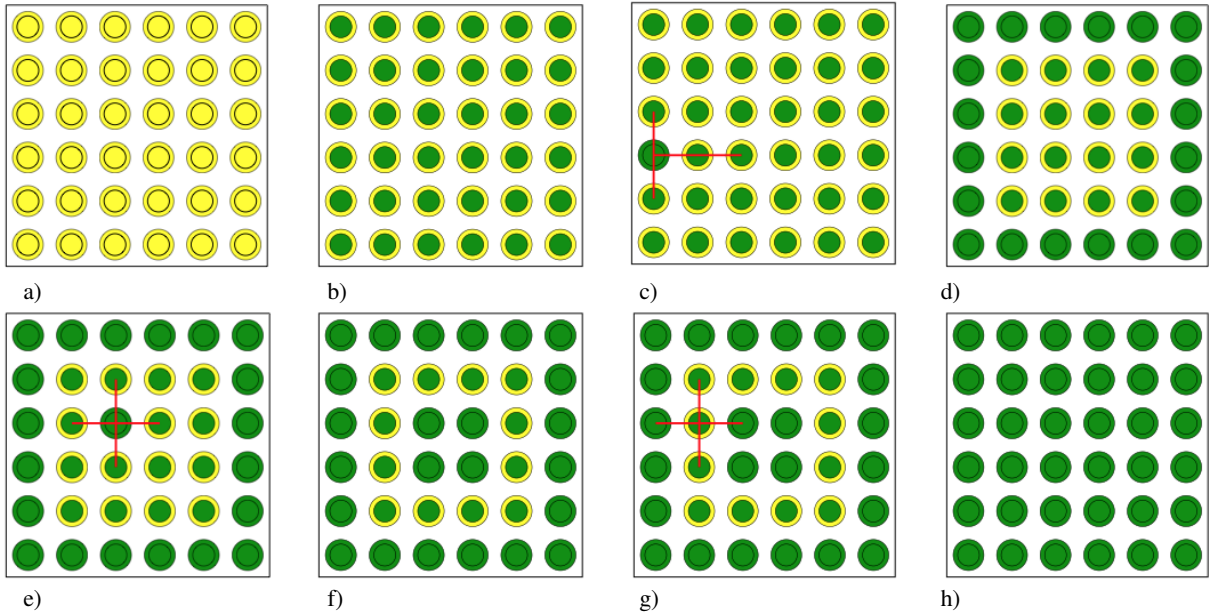

After the acceleration at time $t$ is known at all lattice points, the displacement as well as the distribution functions for the next time step need to be computed in a consistent manner. In this context, we regard the distribution functions, and consequently $\psi$ and $\phi$, to be consistent with the updated displacement field if 
\begin{align}
\begin{split}
    \sum_{\alpha=0}^4 f_{\psi}^{\alpha}(\vec{x},t+\Delta{t}) = \psi(\vec{x},t+\Delta{t})&\\
    \approx (\nabla\times\vec{u})&\vert_{(\vec{x},t+\Delta{t})}, \\
    \sum_{\alpha=0}^4 f_{\phi}^{\alpha}(\vec{x},t+\Delta{t}) = \phi(\vec{x},t+\Delta{t})&\\
    \approx (\nabla\cdot{\vec{u}})&\vert_{(\vec{x},t+\Delta{t})}.
    \label{eq:consistency}
\end{split}
\end{align}
Herein, $(*)\vert_{(\vec{x},t+\Delta{t})}$ means that the spatial derivative $*$ is performed at the lattice point $\vec{x}$ and time $t+\Delta{t}$ via second order accurate finite differences, which is a non-local operation that also involves neighbor lattice points to $\vec{x}$. Fig.~\ref{fig:boundary_update} displays the utilized stencils for this operations as red lines.

The starting point for the algorithm is the situation after the previous time step has been completed as depicted for a quadratic domain in Fig.~\ref{fig:boundary_update}~a). In this figure, lattice points are represented as circles. In order to illustrate the strategy of obtaining consistent displacement and distribution functions, the color of the \textsl{inner} circles also represents the state of the displacement field, i.e. the not yet updated state $\vec{u}(*,t)$ is represented by yellow and the updated state  $\vec{u}(*,t+\Delta{t})$ is indicated by green color. Similarly, the color of the \textsl{outer} circle indicates the state of the rotation $\psi$ and dilatation $\phi$. A yellow outer circle indicates that rotation and dilatation are not updated yet, i.e. the state $\left[\psi(*,t), \phi(*,t)\right]$, whereas a green outer circle indicates the updated state $\left[\psi(*,t+\Delta{t}), \phi(*,t+\Delta{t})\right]$.

The acceleration at all lattice points is known from the Navier-Cauchy equation (\ref{eq:navier}) or the boundary conditions (\ref{eq:acc_neumann}) and (\ref{eq:acc_dirichlet}), which allows to update the displacement field  via (\ref{eq:integration}) as a next step. This leads to an inconsistent situation where the displacement is already updated, but the rotation and the dilatation fields are not, see Fig.~\ref{fig:boundary_update}~b).

The rotation and dilatation fields in the interior are updated via the LBM for the wave equation. This works fine in the interior, where we have pointed out that the Navier-Cauchy equation and the wave equations (\ref{eq:wave-eqs}) are equivalent. Consequently the update of the displacement field by (\ref{eq:navier}) and ($\ref{eq:integration}$) on the one hand, and the update of the distribution functions by  (\ref{eq:LBE}) and the derived rotation and the dilatation by and (\ref{eq:interpretation}) on the other hand are consistent within the limits of the LBM by Chopard et al. \cite{chopard_lattice_1998},  see Fig.~\ref{fig:boundary_update}~e) and f).

However, at boundary lattice points the displacement field is updated from the boundary conditions and the distribution functions can only partially be updated via the LBE. 

Moreover, the neighbors of boundary lattice points, the `second row' boundary lattice points, also cannot be in a consistent state in the sense of (\ref{eq:consistency}) since the finite difference approximation of $\nabla\times\vec{u}$ and $\nabla\cdot\vec{u}$ at those points depends on the displacement of boundary lattice points which in turn is determined only by the boundary conditions and not by the Navier-Cauchy equation.

Thus, in order to model the boundary conditions for the LBM correctly, it is necessary to accomplish two things:
\begin{itemize}
    \item Setting the distribution functions at boundary lattice points consistent with (\ref{eq:consistency}).
    \item Modifying the LBE at ‘second row' lattice points in such a way that consistency is achieved at these points in the sense of (\ref{eq:consistency}).
\end{itemize}
The first requirement is satisfied by setting 
\begin{align}
    &\psi(\vec{x}_k,t+\Delta{t})=(\nabla\times\vec{u})\vert_{(\vec{x}_k,t+\Delta{t})}, \nonumber\\ &\phi(\vec{x}_k,t+\Delta{t})=(\nabla\cdot\vec{u})\vert_{(\vec{x}_k,t+\Delta{t})}
    \label{eq:psi_phi_at_bound}
\end{align}
at boundary lattice points $\boldsymbol{x}_k$ , see Fig.~\ref{fig:boundary_update}~c) and d), and 
\begin{align}
\begin{split}
    f_{\psi\vert\phi}^{0}(\vec{x}_k,t+\Delta{t}) = & a_{0,{\psi\vert\phi}} (\psi\vert\phi)(\vec{x}_k,t+\Delta{t}), \\  f_{\psi\vert\phi}^{\alpha}(\vec{x}_k,t+\Delta{t}) = & a_{{\psi\vert\phi}}(\psi\vert\phi)(\vec{x}_k,t+\Delta{t})\\& +  b \dfrac{\vec{c}^\alpha\cdot\vec{J}_{\psi\vert\phi}(\vec{x}_k,t)}{2c^2}\\
    &\text{for } \alpha\neq{0} .
    \label{eq:correction-f}
\end{split}
\end{align}

In order to fulfill the second requirement, we envision that all changes of $\psi$ and $\phi$ at a boundary lattice points over one time step $t\to{t+\Delta{t}}$ are transported as waves to its neighbors. Assuming that a linear change in time is a reasonable approximation, the average state of a boundary lattice point during this transition is given by the average of its distribution functions at two discrete time steps, i.e.
\begin{equation}
    \tilde{f}^\alpha(\vec{x}_k,\tilde{t}) = \dfrac{1}{2}\left(f_{\psi\vert\phi}^{\alpha}(\vec{x}_k,t+\Delta{t}) + f_{\psi\vert\phi}^{\alpha}(\vec{x}_k,t)\right).
\end{equation}
The transport of this intermediate state to the boundary conditions is obtained by the modified LBE 
\begin{align}
    f^{\alpha}_{\psi\vert\phi} \left( \vec{x} + \vec{c}^{\alpha} \Delta{t}, t + \Delta{t} \right) & = \nonumber\\
    \hat{f}^{\alpha}_{\psi\vert\phi} \argum - \frac{1}{\tau} &\left[ \hat{f}^{\alpha}_{\psi\vert\phi} \argum - f^{\alpha}_{\text{eq},\psi\vert\phi} \argum \right],
    \label{eq:LBE_mod}
\end{align}
where 
\[
    \hat{f}^{\alpha}_{\psi\vert\phi}= 
\begin{cases}
    \tilde{f}^{\alpha}_{\psi\vert\phi}(\vec{x},\tilde{t}),& \text{if } \vec{x} \text{ is boundary lattice point}\\
    {f}^{\alpha}_{\psi\vert\phi}(\vec{x},{t}),              & \text{otherwise.}
\end{cases}
\]
Thus, the modification is only employed if  $\vec{x} + \vec{c}^{\alpha} \Delta{t}$ is a `second row' boundary lattice point. This is not exact, but it is an approximation that leads to a reasonable state of `second row' boundary lattice points, see Fig.~\ref{fig:boundary_update}~g) and h). 

\subsection{Periodic Synchronization}
\label{sec:synchronization}

The proposed LBM for plane strain is  susceptible to instabilities if the dilatation and rotation fields $\psi$ and $\phi$ become inconsistent with the displacements $\boldsymbol{u}$. Since there is no inherent synchronization of these fields, rounding errors are amplified over time and eventually the computed acceleration is sufficiently misaligned with the actual displacement such that the Navier-Cauchy equation (\ref{eq:navier}) is violated. Small inconsistencies originate from the handling of the boundary conditions as described in the last paragraphs of the previous section. 

In order to remedy this problem, we introduce another step in the LBM algorithm, that periodically (every $l$-th timestep where $l\gg{1}$) computes $\psi$ and $\phi$ directly from the displacement field with a finite difference approximation of  (\ref{eq:phi-psi}). As soon as $\psi(\vec{x}_k,t_l)$ and $\phi(\vec{x}_k,t_l)$ are known,
the distribution functions are also corrected according to~(\ref{eq:correction-f}) and the algorithm continues normally in the next time step.

\section{Numerical Examples}

\begin{figure}[htb]
    \centering
    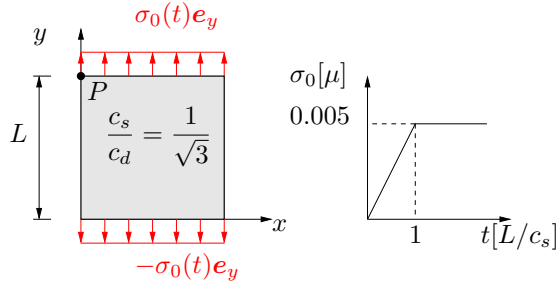
    \caption{A square domain subjected to a tensile load. The right plot displays the applied stress $\sigma_{0}(t)$ as a function of time.}
    \label{fig:tension_setup}
\end{figure}

\begin{figure}[htb]
    \centering
    \includegraphics[width=0.5\textwidth]{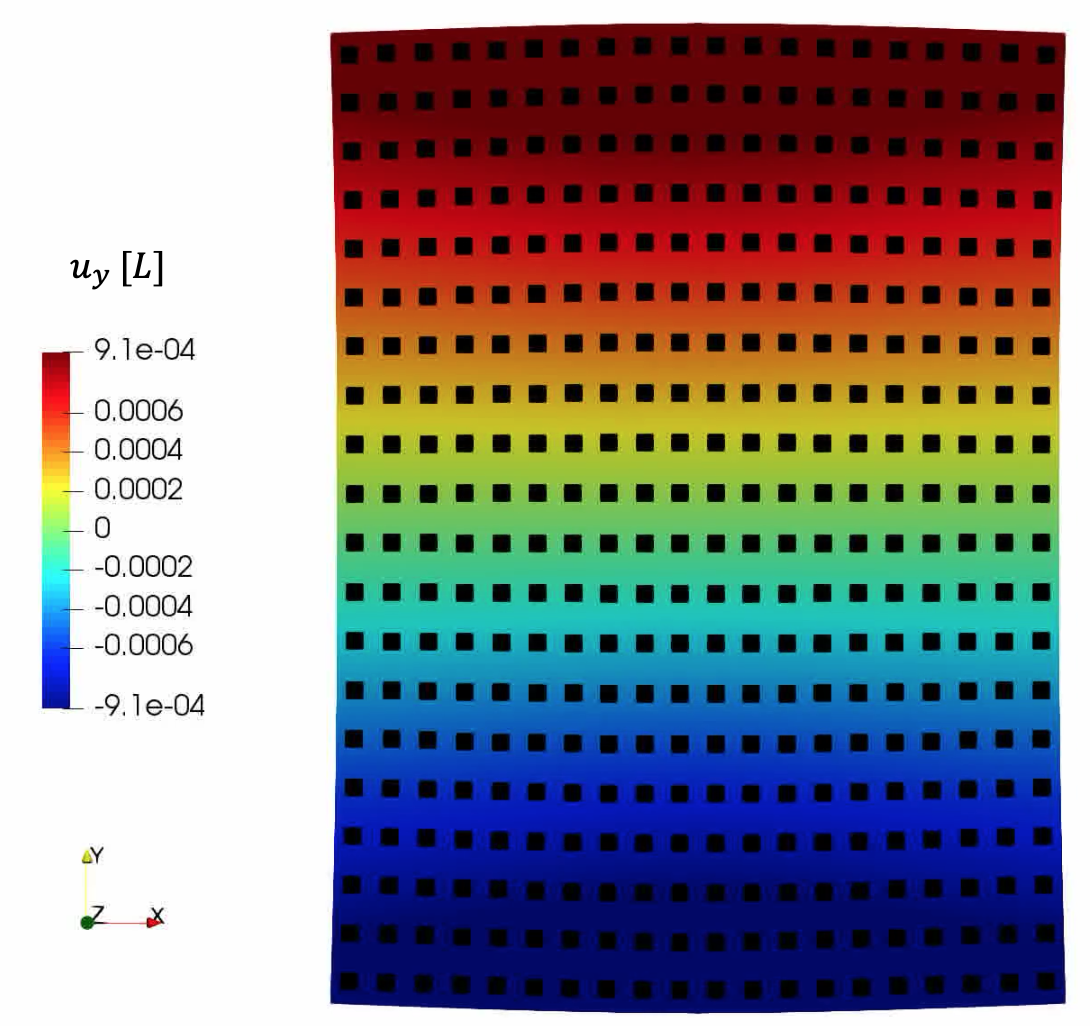}
    \caption{Deformed heat map for a square domain subjected to a tensile load  at time $t=\nicefrac{L}{c_s}$. The deformation is scaled by factor 100. The heat map displays the FEM benchmark results, whereas the black squares indicate the displaced positions of the lattice points.}
    \label{fig:tension_results_1}
\end{figure}

\begin{figure*}[htb]
    \centering
    \includegraphics[width=0.8\textwidth]{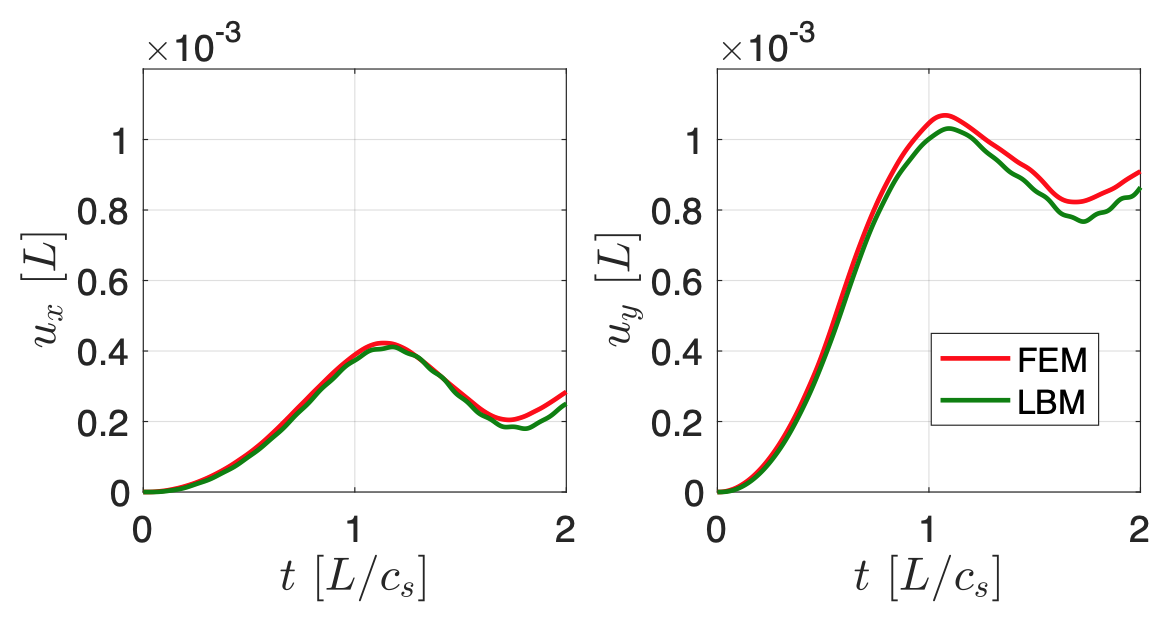}
    \caption{Displacement at the top left corner $P$ of a square domain subjected to a tensile load.}
    \label{fig:tension_results_2}
\end{figure*}

\begin{figure}[htb]
    \centering
    \includegraphics[width=0.4\textwidth]{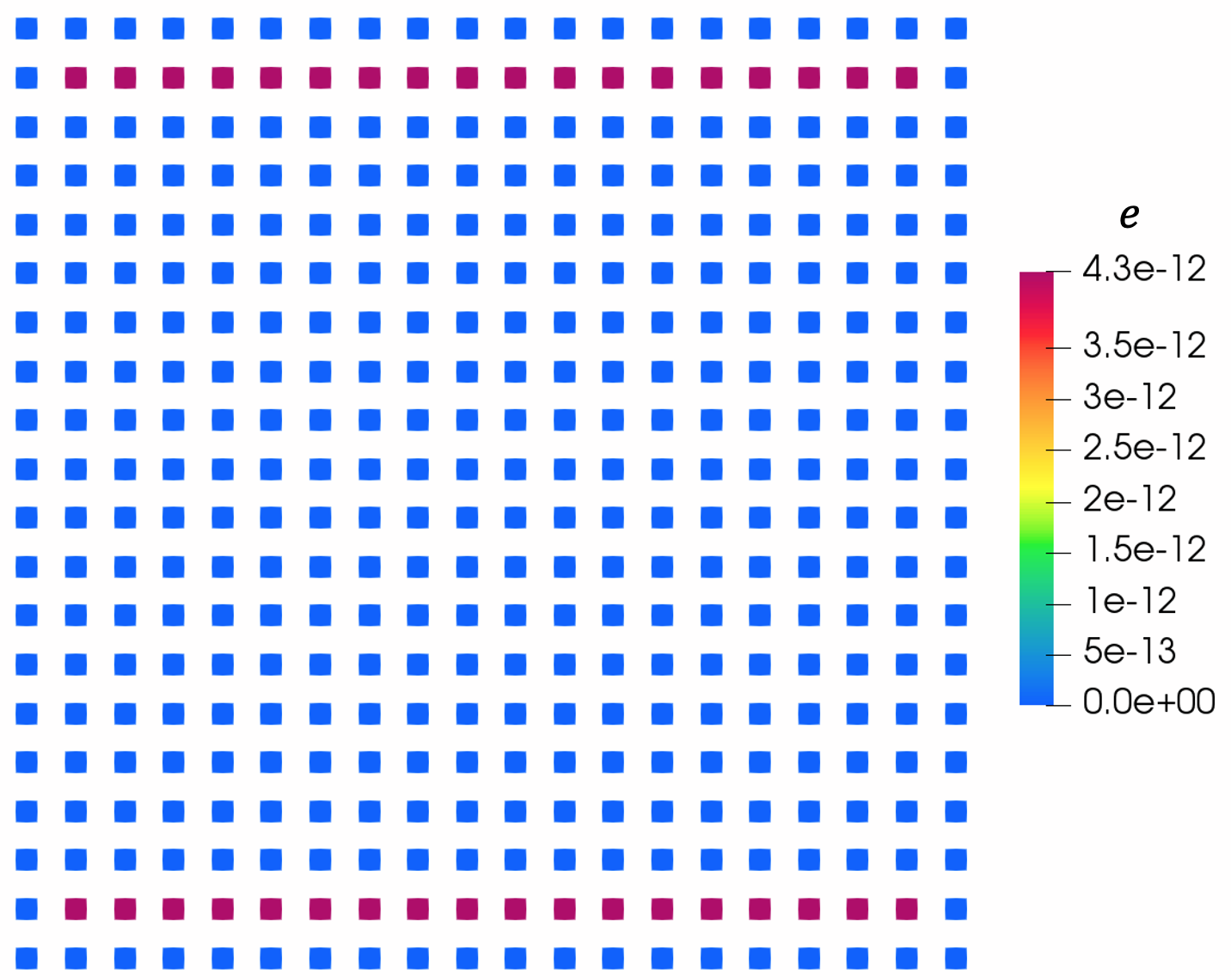}
    \caption{Error $e$ for a square domain subjected to a tensile load at $t\approx0.002\nicefrac{L}{c_s}$. The maximum error is $4.3 \cdot 10^{-12}$ which coincides with the maximum value displayed in the legend color scheme, i.e. dark red.}
    \label{fig:error}
\end{figure}

In order to demonstrate the performance of the proposed LBM, we perform several numerical experiments in which the LBM is compared to results obtained via the established Finite Element Method (FEM). The experiments also demonstrate that the proposed LBM successfully solves boundary value problems that are prevalent in engineering practice and is not restricted to  often rather academic types of boundary value problems, e.g. with periodic boundary conditions.

In all numerical examples, we formulate the problem in terms of  the ratio of wave speeds $\nicefrac{c_s}{c_d}=\nicefrac{1}{\sqrt{3}}$, the wave speed $c_s$, the shear modulus $\mu$, the length scale $L$, and the reference displacement which is also set to $L$. The parameters of the equilibrium distribution functions are defined by (\ref{eq:feq_req}) and (\ref{eq:conditions-parameters}), where a rather extreme value of $a_{0,\phi}=0.9999$ has been found to be required for sufficient stability. Note that setting $a_{0,\phi}=0.9999$ also severely reduces the time step.

The benchmark FEM simulations are performed with bi-linear finite elements and implicit time integration via the standard Newmark method.

\subsection{Tension}

For the first numerical example, a square domain is subjected to a time-dependent, tensile traction $\boldsymbol{t}^*=\pm{\sigma}_0\boldsymbol{e}_y$ load at the top and bottom edges, see Fig.~\ref{fig:tension_setup}. The load is linearly increased from $\sigma_0(t=0)=0$ to \mbox{$\sigma_0(t=\nicefrac{L}{c_s})=0.005\mu$} and held constant afterwards. In this simulation, no periodic synchronization, see section~\ref{sec:synchronization}, is employed. In order to study the performance of the LBM algorithm, we compare the LBM results to an FEM simulation. Fig.~\ref{fig:tension_results_1} shows a deformed heat map of the FEM results in the background, whereas black squares indicate the displaced position\footnote{The displaced position of a lattice point $\vec{x}_k$ is only a result of post-processing, i.e. the lattice remains unchanged. It is computed as $\tilde{\vec{x}}_k(t) = s\vec{u}(\vec{x}_k,t)+\vec{x}_k,$ where $s=100$ is the scaling factor.} of the lattice points. Both simulations are evaluated at time $t=\nicefrac{L}{c_s}$ and the deformation is scaled by a factor of 100. 
It can be observed that the LBM matches the FEM results well and predicts phenomena such as lateral contraction accurately. Fig.~\ref{fig:tension_results_2} explicitly displays the displacement of the top left corner $P$, see also Fig~\ref{fig:hole_setup}, and confirms these findings. We can observe an expected dynamic overshoot and low frequency oscillations in both displacement components, that both the FEM and LBM simulations predict. Nonetheless, Fig.~\ref{fig:tension_results_2} also reveals that erroneous higher frequency oscillations occur in the later stages of the LBM simulation, see the green graph in the plot of $u_y$ for $t>1.5\nicefrac{L}{c_s}$, which indicate that a periodic synchronization may be useful.

As discussed above, we assume that inconsistencies in the sense of violations of (\ref{eq:consistency}) are the cause for these instabilities and that they occur primarily at the `second row' boundary points, see also Fig.~\ref{fig:boundary_update}~g). In order to test this hypothesis, an error measure that is in line with (\ref{eq:consistency}) is defined as
\begin{equation}
    e(\boldsymbol{x},t)= \norm{\begin{pmatrix}
     \sum_{\alpha=0}^4 f_{\psi}^{\alpha}(\vec{x},t) - (\nabla\times\vec{u})\vert_{(\vec{x},t)}\\
    \sum_{\alpha=0}^4 f_{\phi}^{\alpha}(\vec{x},t) - (\nabla\cdot{\vec{u}})\vert_{(\vec{x},t)}
    \end{pmatrix}}_2.
\end{equation}
Fig.~\ref{fig:error} displays this error at time $t\approx0.002\nicefrac{L}{c_s}$ after the corresponding time step has been completely processed. It can be observed that inconsistencies indeed occur at the `second row' boundary points at the top and bottom edges. Although the error is small, without periodic synchronization, it amplifies and eventually manifests as oscillations that can be observed in Fig.~\ref{fig:tension_results_2}.

\subsection{Simple Shear}

\begin{figure*}[htb]
    \centering
     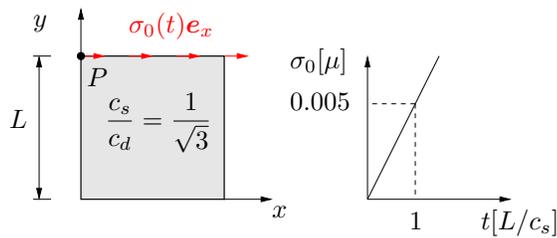
    \caption{A square domain subjected to a shear load. The right plot displays the applied stress $\sigma_{0}(t)$ as a function of time.}
    \label{fig:shear_setup}
\end{figure*}

\begin{figure}[htb]
    \centering
    \includegraphics[width=0.5\textwidth]{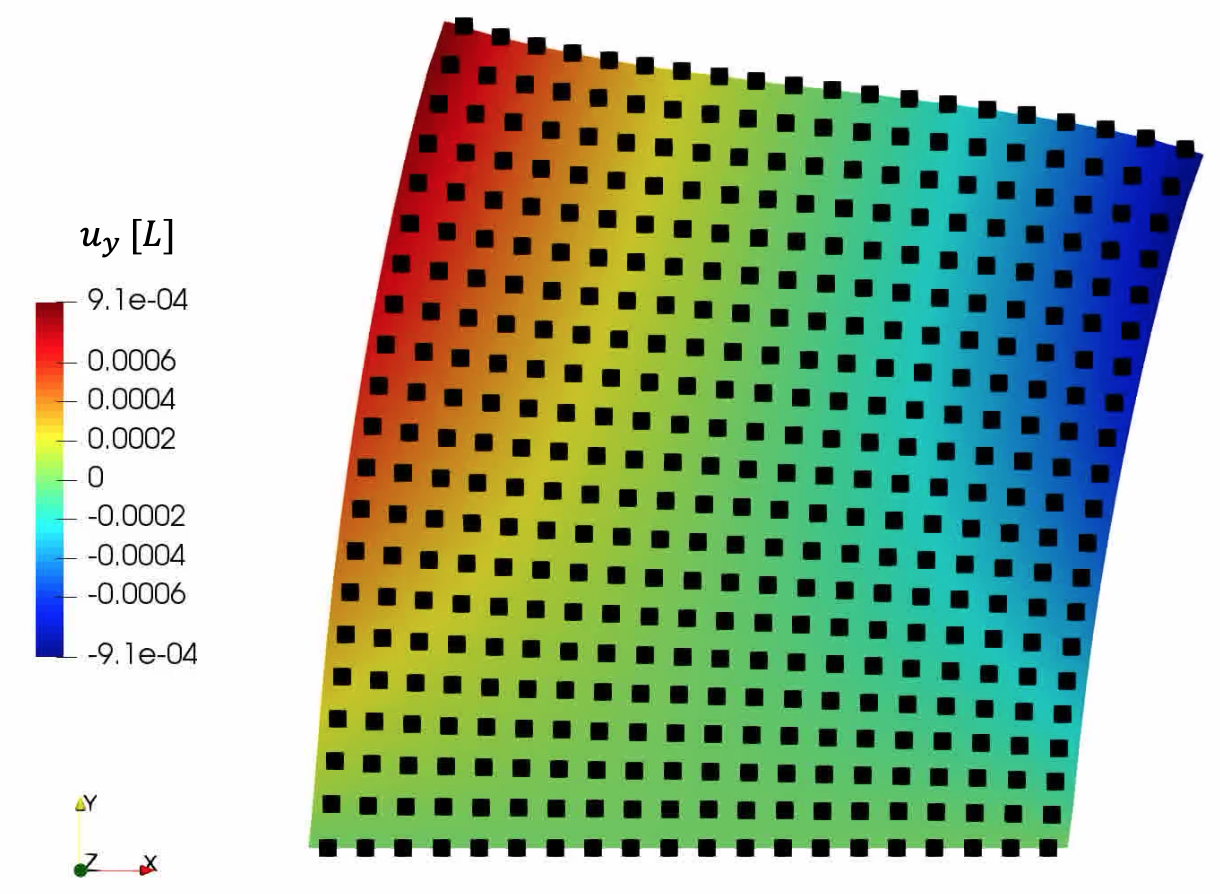}
    \caption{Deformed heat map for a square domain subjected to a shear load  at time $t=\nicefrac{L}{c_s}$. The deformation is scaled by factor 100. The heat map displays the FEM benchmark results, whereas the black squares indicate the displaced positions of the lattice points. For the LBM results a periodic synchronization was performed every 50th time step.}
    \label{fig:shear_results_1}
\end{figure}

\begin{figure*}[htb]
    \centering
    \includegraphics[width=0.8\textwidth]{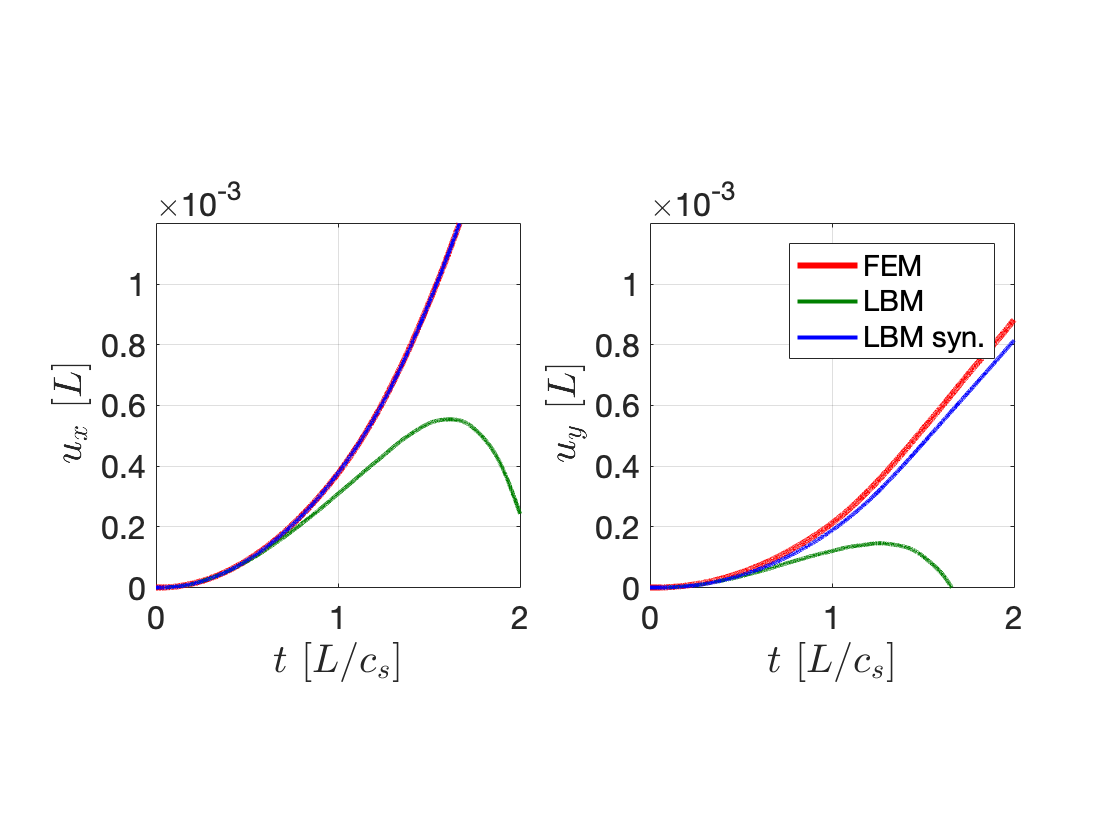}
    \caption{Displacement at the top left corner $P$ of a square domain subjected to a shear load.}
    \label{fig:shear_results_2}
\end{figure*}

The second numerical example uses the same geometric configuration, but differs in terms of the applied boundary conditions. The top edge is subjected to a shear traction that is linearly increased over time, i.e. $\sigma_0=0.005t\mu\nicefrac{c_s}{L}$, see Fig.~\ref{fig:shear_setup}. The bottom edge is subjected to homogeneous Dirichlet boundary conditions $w(x,y=\nicefrac{L}{2},t)=0$. Furthermore, the LBM simulations are run with a periodic synchronization every $50^{\rm th}$ time step and without synchronization. Fig.~\ref{fig:shear_results_1}  displays a deformed (scaled by factor 100) heat map of the FEM results in the background and the displaced lattice points as black squares of the synchronized LBM simulation in the foreground at time $t=\nicefrac{L}{c_s}$. The LBM accurately captures the shear deformation as well. However, as can be observed in Fig.~\ref{fig:shear_results_2}, in this experiment it is strictly necessary to employ the synchronization step, since the LBM simulations without synchronization differ severely from the FEM benchmark after $t\approx0.7\nicefrac{L}{c_s}$.

\subsection{Plate with a Circular Hole}

\begin{figure*}[htb]
    \centering
    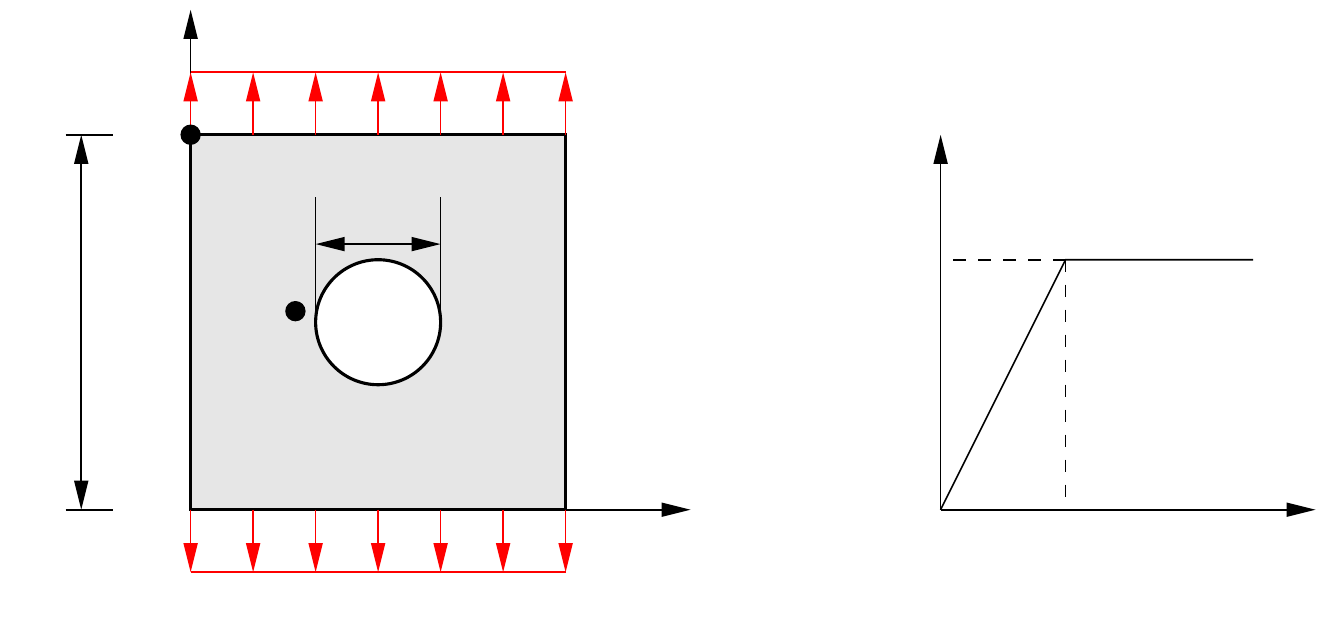
    \caption{A square domain with a hole subjected to a tensile load. Point $Q$ is located at ($-0.175L,0.025L$) relative to a coordinate system which has its origin in the center of the hole. The right plot displays the applied stress $\sigma_{0}(t)$ as a function of time.}
    \label{fig:hole_setup}
\end{figure*}

\begin{figure}[htb]
    \centering
    \includegraphics[width=0.5\textwidth]{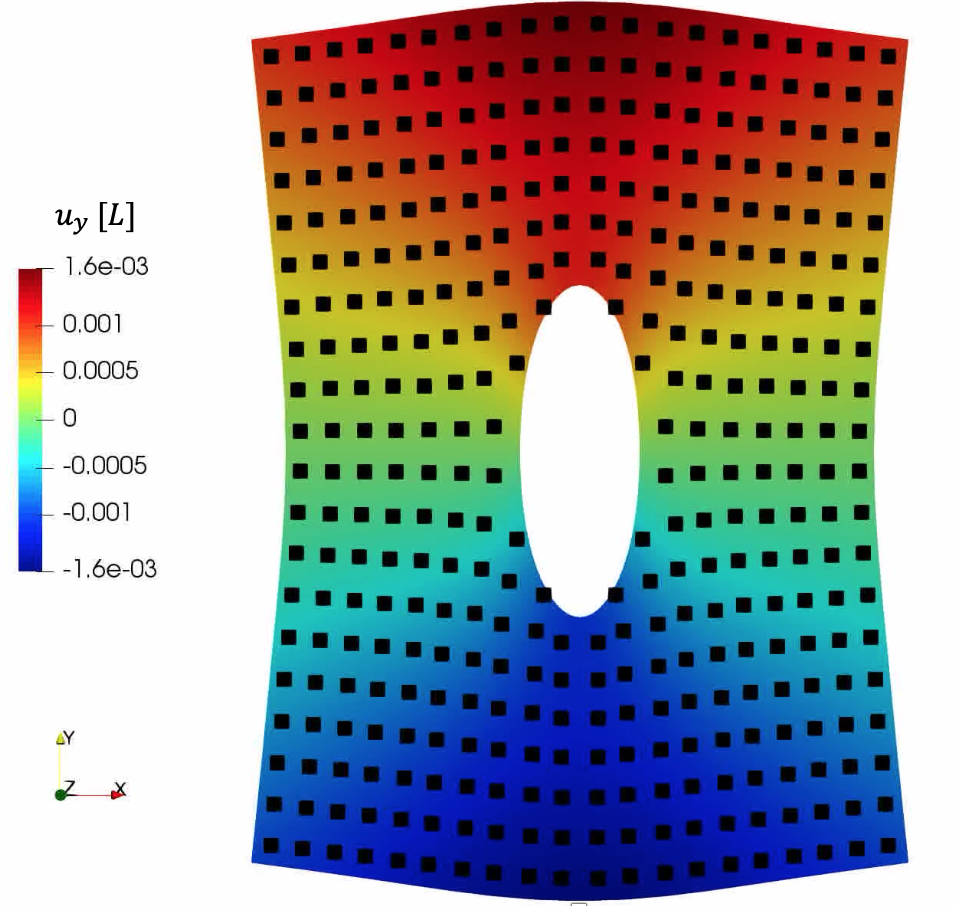}
    \caption{Deformed heat map for a square domain with a hole subjected to a tensile load  at time $t=\nicefrac{L}{c_s}$. The deformation is scaled by factor 100. The heat map displays the FEM benchmark results, whereas the black squares indicate the displaced positions of the lattice points. For the LBM results a periodic synchronization was performed every $50^{\rm th}$ time step.}
    \label{fig:hole_results_1}
\end{figure}

\begin{figure*}[hbtp]
    \centering
    \includegraphics[width=0.8\textwidth]{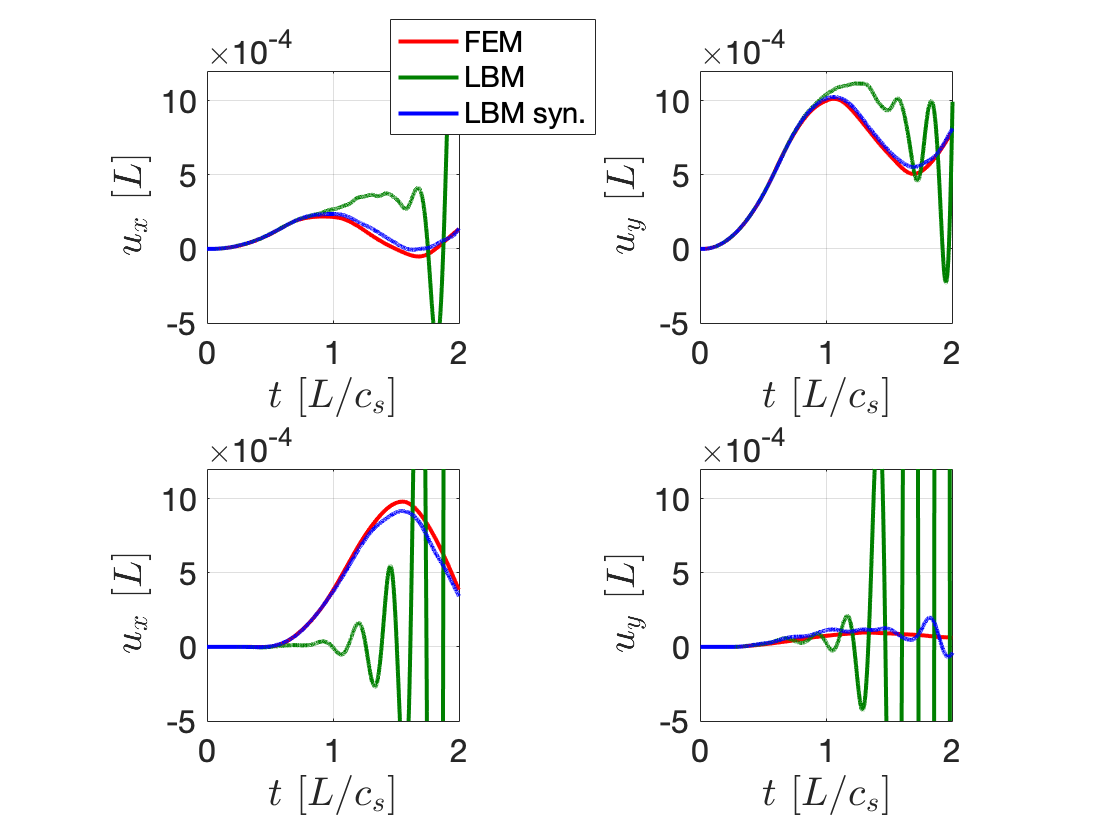}
    \begin{picture}(0,0)
        \put(-360,60){\rotatebox{90}{point $Q$}}
        \put(-360,190){\rotatebox{90}{point $P$}}
    \end{picture}
    \caption{Displacement at the top left corner $P$ (top row) and close to the hole $Q$ (bottom row) of a square domain with a hole subjected to a tensile load.}
    \label{fig:hole_results_2}
\end{figure*}

The last numerical example again considers a square domain that is subjected to a tensile traction load. In order to illustrate the LBMs capabilities to handle non-lattice conforming geometries, the domain includes a circular hole of diameter $0.266L$, see Fig.~\ref{fig:hole_setup}. As in the previous example, we run LBM simulations with periodic synchronization every $50^{\rm th}$ timestep and without any periodic synchronization. Fig~\ref{fig:hole_results_1} displays the scaled deformed configuration for the FEM in the background, as well as for the LBM with synchronization as the black squares in the foreground at time $t=\nicefrac{L}{c_s}$. Again the LBM agrees well with the FEM reference. However, this is only the case if the synchronization is utilized, see Fig.~\ref{fig:hole_results_2}, as the simulation becomes unstable quickly if synchronization is omitted. The oscillations originate from the non-lattice conforming boundaries at the hole as can also be observed in Fig.~\ref{fig:hole_results_2}: the displacement field close to the hole at point $Q$ becomes unstable long before oscillations can be observed at point $P$.

\section{Conclusion}

In this work, a new Lattice Boltzmann Method (LBM) for solving the general plane strain problems is proposed. The plane strain problem is governed by the Navier-Cauchy equation which can be decomposed into two wave equations with different wave speeds for the rotational part of the displacement field and the dilatational part respectively. Based on this observation the new LBM is constructed by employing the established LBM by Chopard et al. \cite{chopard_lattice_1998} to solve the two wave equations separately. 
Chopard et al.'s approach allows enough flexibility to choose the simulated macroscopic wave speed rather independently of the time step and lattice spacing. Thus, 
the proposed method solves both wave equations on the same D2Q5 lattice with the same time discretization. However, this also limits the maximum time step and reduces the computational efficiency of the approach in situations in which larger time step sizes may be feasible.  
The displacement field is eventually obtained by integrating the Navier-Cauchy equation and making use of rotation and dilatation fields computed by the LBM. 

In order to apply Dirichlet and Neumann boundary conditions, a consistent acceleration is computed at boundary lattice points. This is then used for the integration step mentioned above at these points. In order to reconcile the displacements obtained in this way with the LBM quantities such as the rotation and dilatation as well as the distribution functions, the rotation and dilatation fields are computed from a finite difference approximation of the gradient of the displacement field at boundary lattice points. Afterwards, the distribution functions are computed consistently with these rotation and dilatation fields at the boundary points. We mention some of the remaining causes of inconsistencies between rotation and dilatation fields on the one side and the displacement field on the other side. 

These inconsistencies manifest as instabilities in the performed simulations. We address this issue by performing a periodic synchronization in which we compute the rotation and dilatation fields from a a finite difference approximation of the gradient of displacement and subsequently set the distribution functions accordingly.

Lastly, several numerical benchmarks highlight the performance of the new method compared to benchmark FEM simulations. The simulation of a square domain without periodic synchronization under tensile loading shows that the LBM accurately captures simple loading and domains without the synchronization step. However, this example also reveals that the inconsistencies mentioned above indeed occur. 
The second numerical example studies a square domain under simple shear loading conditions. Here, the results only accurately match the FEM simulations if the synchronization step is employed every $50^{\rm th}$ time step. The third numerical example considers the square domain with a hole under tensile load and illustrates that the developed LBM is indeed capable of solving problems in which the geometry does not conform with the lattice, i.e. the boundary does not exactly match the lattice point positions.

We find the performance of the LBM in relation to the FEM promising. However, the periodic synchronization step as well as the rather fine time discretization, that is dictated by the method remains unsatisfactory. In future work, we envision to investigate alternative LBM approaches, but we also want to address the shortcomings of the present LBM by refining the treatment of boundary conditions and exploring the possibility of giving up the same time discretization for both simulated wave equations. This would allow us to use larger time steps and thus increase computational efficiency, but this approach will also involve an additional interpolation step between time steps.

\section*{Acknowledgments}
Open access funding enabled and organized by Projekt DEAL.
The authors gratefully acknowledge the funding by the German Research Foundation (DFG) within the project 423809639.

\bibliographystyle{unsrt}  

\end{document}